%% file: main.tex
\documentclass[12,ptoneside]{amsart}

\usepackage{graphicx}        % standard LaTeX graphics tool
\usepackage{multicol}        % used for the two-column index
\usepackage[bottom]{footmisc}% places footnotes at page bottom

\usepackage{url}
\usepackage{subfig}

\usepackage{amsmath,amstext,amssymb,bbm,amsbsy}
\usepackage{mathtools}

\newcommand{\be}{\begin{equation}}
\newcommand{\ee}{\end{equation}}
\newcommand{\benn}{\begin{equation*}}
\newcommand{\eenn}{\end{equation*}}

\newtheorem{proposition}{Proposition.}%[section]

\DeclareMathOperator*{\argmin}{arg\,min}

\newcommand{\diff}{\mathop{}\!d}

\newcommand{\normsq}[1]{\left\Vert{#1} \right\Vert_{2}^{2}}

\def\boldeta {\boldsymbol{\eta}}
\def\bx {\boldsymbol{x}}
\def\by {\boldsymbol{y}}
\def\bv {\boldsymbol{v}}
\def\bu {\boldsymbol{u}}

\def\bp {\boldsymbol{p}}
\def\R {\mathbb{R}}

\def\dom {\mathrm{dom~}}
\def\inter {\mathrm{int~}}

\def\Rn {\R^n}
\def\gmRn {\Gamma_0(\R^n)}
\def\Laplacian{\Delta_{\boldsymbol x}}

\date{\today. This work was funded by NSF 1820821. Authors' names are given in last/family name alphabetical order.}

\begin{document}

\title{Connecting Hamilton--Jacobi partial differential equations with maximum a posteriori and posterior mean estimators for some non-convex priors}
\author[J. Darbon]{J\'er\^ome Darbon}
\address{Division of Applied Mathematics, Brown University.}
\email{jerome\_darbon@brown.edu}

\author[G. P. Langlois]{Gabriel P. Langlois}
\address{Division of Applied Mathematics, Brown University.}
\email{gabriel\_provencher\_langlois@brown.edu}

\author[T. Meng]{Tingwei Meng}
\address{Division of Applied Mathematics, Brown University.}
\email{tingwei\_meng@brown.edu}

%\author{J\'er\^ome Darbon and Gabriel P. Langlois and Tingwei Meng}
%\institute{Authors are ordered by alphabetical family names\\J\'er\^ome Darbon \at Division of Applied Mathematics, Brown University, \email{jerome_darbon@brown.edu}
%\and Gabriel P. Langlois \at  Division of Applied Mathematics, Brown University \email{gabriel_provencher_langlois@brown.edu}
%\and Tingwei Meng \at  Division of Applied Mathematics, Brown University \email{tingwei_meng@brown.edu}
%}
%This work was funded by NSF 1820821.
\maketitle

\begin{abstract}
Many imaging problems can be formulated as inverse problems expressed as finite-dimensional optimization problems. These optimization problems generally consist of minimizing the sum of a data fidelity and regularization terms. In~\cite{darbon2015convex,darbon2019decomposition}, connections between these optimization problems and (multi-time) Hamilton--Jacobi partial differential equations have been proposed under the convexity assumptions of both the data fidelity and regularization terms. In particular, under these convexity assumptions, some representation formulas for a minimizer can be obtained. From a Bayesian perspective, such a minimizer can be seen as a maximum a posteriori estimator.
In this chapter, we consider a certain class of non-convex regularizations and show that similar representation formulas for the minimizer can also be obtained. This is achieved by leveraging min-plus algebra techniques that have been originally developed for solving certain Hamilton--Jacobi partial differential equations arising in optimal control. Note that connections between viscous Hamilton--Jacobi partial differential equations and Bayesian posterior mean estimators with Gaussian data fidelity terms and log-concave priors have been highlighted in~\cite{darbon2020bayesian}. We also present similar results for certain Bayesian posterior mean estimators with Gaussian data fidelity and certain non-log-concave priors using an analogue of min-plus algebra techniques.
\end{abstract}

\input{introduction}
\input{firstorder}

\input{secondorder}
\input{conclusion}

\bibliographystyle{siam}
\bibliography{biblio,biblist}
\end{document}

%% file: introduction.tex
\section{Introduction}
Many low-level signal, image processing and computer vision problems are formulated as inverse problems that can be solved using variational \cite{aubert.02.book,2009Scherzer,veseGuyager.16.book} or Bayesian approaches~\cite{winkler.03.book}. Both approaches have been very effective, for example, at solving image restoration~\cite{bouman.93.itip,likasetal.04.ieeesp,rudin.92.phys}, segmentation~\cite{boykovetal.01.pami,chanetal.06.siam,chanvese.01.itip} and image decomposition problems \cite{aujol2005image,OSV}.

As an illustration, let us consider the following image denoising problem in finite dimension that formally reads as follows:
$$
\bx = \bar{\bu} + \boldeta,
$$
where $\bx \in \Rn$ is the observed image that is the sum of an unknown ideal image $\bar{\bu} \in \Rn$ and an additive perturbation or noise realization $\boldeta \in \Rn$. We aim to estimate~$\bar{\bu}$.

A standard variational approach for solving this problem consists of estimating $\bar{\bu}$ as a minimizer of the following optimization problem 
\begin{equation}
\label{eq:canonical-optimization-pb}
\min_{\bu \in \Rn} \left\{\lambda D(\bx-\bu) + J(\bu)\right\},
\end{equation}
where $D:\Rn \to \R$ is generally called the data fidelity term and contains the knowledge we have on the perturbation $\boldeta$, while $J:\Rn \to \R \cup \{+\infty\}$ is called the regularization term and encodes the knowledge on the image we wish to reconstruct. The non-negative parameter $\lambda$ relatively weights the data fidelity and the regularization terms. Note that minimizers of~\eqref{eq:canonical-optimization-pb} are called maximum a posteriori (MAP) estimators in a Bayesian setting. Also note that variational-based approaches for estimating $\bar{\bu}$ are particularly appealing when 
both the data fidelity and regularization terms are convex because~\eqref{eq:canonical-optimization-pb} becomes a convex optimization problem that can be efficiently solved using convex optimization algorithms (see e.g.,~\cite{chambolle2016introduction}). 
Many regularization terms have been proposed in the literature~\cite{aubert.02.book,winkler.03.book}. Popular choices for these regularization terms involve robust edge preserving priors
\cite{bouman.93.itip,charbonnieretal.97.ip,geman_yang.95.ip,gemanReynolds.92.pami,NikolovaChan-07-ip,nikolva-05-sc,rudin.92.phys} because they allow the reconstructed image to have sharp edges. For the sake of simplicity, we only describe in this introduction regularizations that are expressed using pairwise interactions which take the following form
\begin{equation}
\label{eq:pairwiseInteractions}
J(\bu) = \sum_{i,j=1}^n w_{ij}  f(u_i - u_j),
\end{equation}
where $f : \R \to \R \cup \{+\infty\}$ and $w_{i,j} \geqslant 0$. Note that our results that will be presented later do not rely on pairwise interaction-based models and work for more general regularization terms. A popular choice is the celebrated Total Variation~\cite{bouman.93.itip,rudin.92.phys}, which corresponds to consider $f(z) = |z|$ in \eqref{eq:pairwiseInteractions}. The use of Total Variation as a regularization term has been very popular since the seminal works of~\cite{bouman.93.itip,rudin.92.phys} because it is convex and allows the reconstructed image to preserve edges well. When the data fidelity $D$ is quadratic, this model is known as the celebrated Rudin-Osher-Fatemi model~\cite{rudin.92.phys}. Following the seminal works of~\cite{charbonnieretal.97.ip,geman_yang.95.ip,gemanReynolds.92.pami}, another class of edges preserving priors corresponds to half-quadratic-based regularizations that read as follows:
\begin{equation}\label{eqt:half_quadratic}
f(z) = 
\begin{dcases}
|z|^2 & \mbox{ if } |z|\leqslant 1,\\
1 & \mbox{ otherwise}.
\end{dcases}
\end{equation}
Note that the quadratic term above can be replaced by $|\cdot|$, i.e., we consider
\begin{equation}\label{eqt:l1}
f(z) 
= \begin{dcases}
|z| & \mbox{ if } |z| \leqslant 1,\\
1 & \mbox{otherwise},
\end{dcases}
\end{equation}
which corresponds to the truncated Total Variation regularization (see~\cite{darbonetal.09.icip,Douetal.17.ieee} for instance).

There is a large body of literature on variational methods (e.g.,~\cite{aubert.02.book,Chambolle10anintroduction,chan2005image,2009Scherzer,veseGuyager.16.book}). In particular,
in \cite{darbon2015convex,darbon2019decomposition}, connections between convex optimization optimization problems of the form of~\eqref{eq:canonical-optimization-pb} and 
Hamilton--Jacobi partial differential equations (HJ PDEs) were highlighted. Specifically, it is shown that the dependence of the minimal value of these problems with respect to the observed data $\bx$ and the smoothing parameter $\lambda$ are governed by HJ PDEs, where the initial data corresponds to the regularization term $J$ and the Hamiltonian is related to the data fidelity (see Section~\ref{sec:firstorder} for details).  However, the connections between HJ PDEs and certain variational imaging problems described in~\cite{darbon2015convex,darbon2019decomposition} require the convexity of both the data fidelity and regularization terms. Note that these connections between HJ PDEs and imaging problems also hold for image decomposition models~(see Section~\ref{subsec:MultiTimeHJ}) using multi-time HJ PDEs~\cite{darbon2019decomposition}.

Our goal is to extend the results of ~\cite{darbon2015convex,darbon2019decomposition} to certain non-convex regularization terms using min-plus algebra techniques~\cite{akian2006max,akian2008max, dower2015max,Fleming2000Max,gaubert2011curse,Kolokoltsov1997Idempotent,mceneaney2006max,McEneaney2007COD,mceneaney2008curse,mceneaney2009convergence} that were originally designed for solving certain HJ PDEs arising in optimal control problems. We also propose an analogue of this approach for certain Bayesian posterior mean estimators when the data fidelity is Gaussian.

The rest of this chapter is as follows. Section~\ref{sec:firstorder} reviews connections of image denoising and decomposition models with HJ PDEs under convexity assumptions. We then present a min-plus algebra approach for single time and multi-time HJ PDEs that allows us to consider certain non-convex regularizations in these image denoising and decomposition models. In particular, this min-plus algebra approach yields practical numerical optimization algorithms for solving certain image denoising and decomposition models.
Section~\ref{sec:viscousHJ} reviews connections between viscous HJ PDEs and posterior mean estimators with Gaussian data fidelity term and log-concave priors. We also present an analogue of the min-plus algebra technique for these viscous HJ PDEs with certain priors that are not log-concave.
Finally, we draw some conclusions in Section~\ref{sec:conclusion}.

%% file: firstorder.tex
\section{First order Hamilton--Jacobi PDEs and optimization problems}
\label{sec:firstorder}

In this section, we discuss the connections between some variational optimization models in imaging sciences and HJ PDEs. In subsection~\ref{subsec:singleTimeHJ}, we consider the convex image denoising model \eqref{eq:canonical-optimization-pb} and the single time HJ PDE. In subsection~\ref{subsec:MultiTimeHJ}, we review the connections between convex image decomposition models and the multi-time HJ PDE system. In subsection~\ref{subsec:minplusMultiTimeHJ}, we use the min-plus algebra technique to solve certain optimization problems in which one regularization term is non-convex. In subsection~\ref{sec: eg}, we provide an application of the min-plus algebra technique to certain image decomposition problems, which yields practical numerical optimization algorithms.

\subsection{Single time HJ PDEs and image denoising models}
\label{subsec:singleTimeHJ}

As described in the introduction, an important class of optimization models in imaging sciences for denoising takes the form of~\eqref{eq:canonical-optimization-pb}, where $\lambda>0$ is a positive parameter, $\bx\in\R^n$ is the observed image with $n$ pixels, and $\bu\in \R^n$ is the reconstructed image. The objective function is the weighted sum of the convex regularization term $J$ and the convex data fidelity term $D$. 

The connection between the class of optimization models~\eqref{eq:canonical-optimization-pb} and first order HJ PDEs has been discussed in~\cite{darbon2015convex}. Specifically, if the data fidelity term $\lambda D$ can be written in the form of $t H^*\left(\frac{\cdot}{t}\right)$ (where $H^*$ denotes the Legendre transform of a convex function $H$ and $t>0$ is a new parameter that depends on $\lambda$), then the minimization problem \eqref{eq:canonical-optimization-pb} defines a function $S\colon\Rn \times (0,+\infty) \to \R$ as follows
\begin{equation} \label{eqt: 1order_LO}
S(\bx,t) = \min_{\bu\in\R^n} \left\{J(\bu) + t H^*\left(\frac{\bx - \bu}{t}\right)\right\}.
\end{equation}
For instance, if the noise is assumed to be Gaussian, independent, identically distributed and additive, we impose the quadratic data fidelity $D(\bx) = \frac{1}{2}\|\bx\|_2^2$ for each $\bx\in\R^n$. Then $D$ satisfies $\lambda D(\bx) = tH^*\left(\frac{\bx}{t}\right)$ where $H^*(\bx) = \frac{1}{2}\|\bx\|_2^2$ and $t = \frac{1}{\lambda}$.

Formula \eqref{eqt: 1order_LO} is called the Lax-Oleinik formula \cite{bardi1984hopf,evans1998partial,hopf1965generalized} in the PDE literature, and it solves the following first order HJ PDE 
\begin{equation}\label{eqt: 1order_HJ}
\begin{dcases}
\frac{\partial S}{\partial t}(\bx,t) + H(\nabla_{\bx} S(\bx,t)) = 0 &\bx\in \Rn, t>0,\\
S(\bx,0) = J(\bx) &\bx\in \Rn,
\end{dcases}
\end{equation}
where the function $H\colon \R^n\to\R$ is called the Hamiltonian, and $J\colon \R^n\to\R\cup\{+\infty\}$ is the initial data. In \cite{darbon2015convex}, a representation formula for the minimizer of~\eqref{eqt: 1order_LO} is given, and we state it in the following proposition. Here and in the remainder of this chapter, we use $\Gamma_0(\R^n)$ to denote the set of convex, proper and lower semicontinuous functions from $\R^n$ to $\R\cup\{+\infty\}$.

\begin{proposition} \label{prop: firstorder_1}
Assume $J\in \Gamma_0(\R^n)$, and assume $H\colon \R^n\to \R$ is a differentiable, strictly convex and 1-coercive function. Then the Lax-Oleinik formula \eqref{eqt: 1order_LO} gives the differentiable and convex solution $S\colon \R^n\times (0,+\infty) \to \R$ to the HJ PDE \eqref{eqt: 1order_HJ}.
Moreover, for each $\bx\in\R^n$ and $t >0$, the minimizer in \eqref{eqt: 1order_LO} exists and is unique, which we denote by $\bu(\bx,t)$, and satisfies
\be \label{eqt: prop1_relation_u_p}
\bu(\bx,t) = \bx - t\nabla H(\nabla_{\bx} S(\bx,t)).
\ee
\end{proposition}

Equation \eqref{eqt: prop1_relation_u_p} in this proposition gives the relation between the minimizer $\bu$ in the the Lax-Oleinik formula \eqref{eqt: 1order_LO} and the spatial gradient of the solution to the HJ PDE \eqref{eqt: 1order_HJ}. In other words, one can compute the minimizer in the corresponding denoising model \eqref{eq:canonical-optimization-pb} using the spatial gradient $\nabla_{\bx} S(\bx, t)$ of the solution, and vice versa. 

There is another set of assumptions for the conclusion of the proposition above to hold. For the details, we refer the reader to \cite{darbon2015convex}.

\subsection{Multi-time HJ PDEs and image decomposition models}
\label{subsec:MultiTimeHJ}

In this subsection, we consider the following image decomposition models:
\be \label{eqt: 1order_decomposition}
\min_{\bu_1, \dots, \bu_N\in \R^n} \left\{ J\left(\bx - \sum_{i=1}^N \bu_i \right) + \sum_{i=1}^N \lambda_i f_i(\bu_i) \right\},
\ee
where $\lambda_1, \dots, \lambda_N$ are positive parameters, $\bx\in\R^n$ is the observed image with $n$ pixels, and $\bu_1, \dots, \bu_N \in \R^n$ correspond to the decomposition of the original image $\bx$.
In \cite{darbon2019decomposition}, the relation between the decomposition model \eqref{eqt: 1order_decomposition} and the multi-time HJ PDE system has been proposed under the convexity assumptions of $J$ and the functions $f_1, \dots, f_N$.

In the decomposition model, an image is assumed to be the summation of $N+1$ components, denoted as $\bu_1, \dots, \bu_N$ and the residual $\bx - \sum_{i=1}^N \bu_i$. The feature of each part $\bu_i$ is characterized by a convex function $f_i$, and the residual $\bx - \sum_{i=1}^N \bu_i$ is characterized by a convex regularization term $J$. If the function $\lambda_i f_i$ can be written in the form of $t_i H_i^*\left(\frac{\cdot}{t_i}\right)$ (where $H_i^*$ denotes the Legendre transform of a convex function $H_i$ and $t_i>0$ is a new parameter which depends on $\lambda_i$) for each $i\in \{1,\dots, N\}$, then the image decomposition model \eqref{eqt: 1order_decomposition}
defines a function $S\colon\Rn \times (0,+\infty)^N \to \R$ as follows
\be \label{eqt: multitime_LO}
S(\bx, t_1,\dots, t_N) = \min_{\bu_1,\dots, \bu_N\in\R^n} \left\{ J\left(\bx - \sum_{i=1}^N \bu_i \right) + \sum_{i=1}^N t_i H_i^*\left(\frac{\bu_i}{t_i}\right)\right\}.
\ee
This formula is called the generalized Lax-Oleinik formula \cite{lions1986hopf,tho2005hopf} which solves the following multi-time HJ PDE system
\begin{equation} \label{eqt: multitime_HJ}
\begin{dcases}
\frac{\partial S(\bx, t_1,\dots,t_N)}{\partial t_1} + H_1(\nabla_{\bx} S(\bx, t_1,\dots,t_N)) = 0 &\bx\in\Rn, t_1,\cdots,t_N>0, \\
\quad\quad\quad\;\; \vdots  \\
\frac{\partial S(\bx, t_1,\dots,t_N)}{\partial t_j} + H_j(\nabla_{\bx} S(\bx, t_1,\dots,t_N)) = 0 &\bx\in\Rn, t_1,\cdots,t_N>0, \\
\quad\quad\quad\;\; \vdots  \\
\frac{\partial S(\bx, t_1,\dots,t_N)}{\partial t_N} + H_N(\nabla_{\bx} S(\bx, t_1,\dots,t_N)) = 0 &\bx\in\Rn, t_1,\cdots,t_N>0, \\
S(\bx,0,\cdots,0) = J(\bx) &\bx\in \Rn,
\end{dcases}
\end{equation}
where $H_1, \dots, H_N\colon \R^n\to \R$ are called Hamiltonians, and $J\colon \R^n\to \R\cup\{+\infty\}$ is the initial data.
Under certain assumptions (see Prop.~\ref{prop: firstorder_2} below), the generalized Lax-Oleinik formula~\eqref{eqt: multitime_LO} gives the solution $S(\bx, t_1, \dots, t_N)$ to the multi-time HJ PDE system \eqref{eqt: multitime_HJ}.
In \cite{darbon2019decomposition}, the relation between the minimizer in \eqref{eqt: multitime_LO} and the spatial gradient $\nabla_{\bx} S(\bx, t_1, \dots, t_N)$ of the solution to the multi-time HJ PDE system \eqref{eqt: multitime_HJ} is studied. This relation is described in the following proposition.

\begin{proposition} \label{prop: firstorder_2}
Assume $J \in \Gamma_0(\R^n)$, and assume $H_j \colon \R^n \to \R$ is a convex and 1-coercive function for each $j \in \{1, \dots,N\}$. Suppose there exists $j\in\{1,\dots,N\}$ such that $H_j$ is strictly convex.
Then the generalized Lax-Oleinik formula \eqref{eqt: multitime_LO} gives the differentiable and convex solution $S\colon \R^n\times (0,+\infty)^N \to \R$ to the multi-time HJ PDE system \eqref{eqt: multitime_HJ}.
Moreover, for each $\bx\in\R^n$ and $t_1, \dots, t_N >0$, the minimizer in \eqref{eqt: multitime_LO} exists. We denote by $(\bu_1(\bx,t_1, \dots, t_N), \dots, \bu_N(\bx,t_1, \dots, t_N))$ any minimizer of the minimization problem in \eqref{eqt: multitime_LO} with parameters $\bx\in \R^n$ and $t_1, \dots, t_N \in (0,+\infty)$. Then, for each $j\in\{1,\dots,N\}$, there holds
\benn
\bu_j(\bx,t_1, \dots, t_N) \in t_j\partial H_j (\nabla_{\bx}S(\bx,t_1, \dots, t_N)),
\eenn
where $\partial H_j$ denotes the subdifferential of $H_j$.

Furthermore, if all the Hamiltonians $H_1,\dots, H_N$ are differentiable, then the
minimizer is unique and satisfies
\be \label{eqt: u_formula_multitime}
\bu_j(\bx,t_1, \dots,t_N) = t_j\nabla H_j(\nabla_{\bx} S(\bx,t_1, \dots,t_N)), 
\ee
for each $j\in\{1,\dots, N\}$.
\end{proposition}

As a result, when the assumptions in the proposition above are satisfied, one can compute the minimizer to the corresponding decomposition model \eqref{eqt: 1order_decomposition} using equation \eqref{eqt: u_formula_multitime} and the spatial gradient $\nabla_{\bx} S(\bx,t_1, \dots,t_N)$ of the solution to the multi-time HJ PDE \eqref{eqt: multitime_HJ}.

\subsection{Min-plus algebra for HJ PDEs and certain non-convex regularizations}
\label{subsec:minplusMultiTimeHJ}

In the previous two subsections, we considered the optimization models \eqref{eq:canonical-optimization-pb} and \eqref{eqt: 1order_decomposition} where each term was assumed to be convex. When $J$ is non-convex, solutions to~\eqref{eqt: 1order_HJ} may not be classical (in the sense that it is not differentiable). It is well-known that the concept of viscosity solutions \cite{Bardi1997Optimal,barles1994solutions, BARRON1984213, crandall1992user, evans1998partial,fleming2006controlled} is generally the appropriate notion of solutions for these HJ PDEs. Note that Lax-Oleinik formulas \eqref{eq:canonical-optimization-pb} and \eqref{eqt: 1order_decomposition} yield viscosity solutions to their respective HJ PDEs \eqref{eqt: 1order_HJ} and \eqref{eqt: multitime_HJ}. However, these Lax-Oleinik formulas result in non-convex optimization problems.

In this subsection, we use the min-plus algebra technique \cite{akian2006max,akian2008max, dower2015max,Fleming2000Max,gaubert2011curse,Kolokoltsov1997Idempotent,mceneaney2006max,McEneaney2007COD,mceneaney2008curse,mceneaney2009convergence} to handle the cases when the term $J$ in \eqref{eq:canonical-optimization-pb} and \eqref{eqt: 1order_decomposition} is assumed to be a non-convex function in the following form
\be \label{eqt: minJi}
J(\bx) = \min_{i \in \{1,\dots,m\}} J_i(\bx) \text{ for every } \bx\in\R^n,
\ee
where $J_i \in \Gamma_0(\R^n)$ for each $i\in\{1,\dots, m\}$.

First, we consider the single time HJ PDE \eqref{eqt: 1order_HJ}.
By min-plus algebra theory, the semi-group of this HJ PDE is linear with respect to the min-plus algebra. In other words, under certain assumptions the solution $S$ to the HJ PDE $\frac{\partial S}{\partial t}(\bx,t) + H(\nabla_{\bx} S(\bx,t)) = 0$ with initial data $J$ is the minimum of the solution $S_i$ to the HJ PDE $\frac{\partial S_i}{\partial t}(\bx,t) + H(\nabla_{\bx} S_i(\bx,t)) = 0$ with initial data $J_i$. Specifically, if the Lax-Oleinik formula \eqref{eqt: 1order_LO} solves the HJ PDE \eqref{eqt: 1order_HJ} for each $i\in\{1,\dots,m\}$ and the minimizer $\bu$ exists (for instance, when $J_i\in \Gamma_0(\R^n)$ for each $i\in\{1,\dots,m\}$, and $H\colon \R^n\to \R$ is a differentiable, strictly convex and 1-coercive function), then we have
\be \label{eqt: minplus_singletime}
\begin{split}
S(\bx, t) &= \min_{\bu\in\R^n} \left\{ J(\bu) + t H^*\left(\frac{\bx - \bu}{t}\right)\right\}\\
&= \min_{\bu\in\R^n} \left\{ \min_{i \in \{1,\dots,m\}} J_i(\bu) + t H^*\left(\frac{\bx - \bu}{t}\right)\right\}\\
&= \min_{\bu\in\R^n} \min_{i \in \{1,\dots,m\}}\left\{  J_i(\bu) + t H^*\left(\frac{\bx - \bu}{t}\right)\right\}\\
&=  \min_{i \in \{1,\dots,m\}}\left\{ \min_{\bu\in\R^n}\left\{ J_i(\bu) + t H^*\left(\frac{\bx - \bu}{t}\right)\right\}\right\}\\
&= \min_{i \in \{1,\dots,m\}} S_i(\bx,t).
\end{split}
\ee
Therefore, the solution $S(\bx,t)$ is given by the pointwise minimum of $S_i(\bx,t)$ for $i\in\{1,\dots,m\}$. Note that the Lax-Oleinik formula~\eqref{eqt: 1order_LO}  yields a convex problem for each $S_i(\bx,t)$ with $i\in\{1,\dots,m\}$. Therefore this approach seems particularly appealing to solve these non-convex optimization problems and associated HJ PDEs. Note that such an approach is embarrassingly parallel since we can solve the initial data $J_i$ for each  $i\in\{1,\dots,m\}$ independently and compute in linear time the pointwise minimum. However, this approach is only feasible if $m$ is not too big. We will see later in this subsection that robust edge preserving priors (e.g., truncated Total Variation or truncated quadratic) can be written in the form of~\eqref{eqt: minJi} but $m$ is exponential in $n$.  

\bigbreak

We can also compute the set of minimizers $\bu(\bx,t)$ as follows. Here, we abuse notation and use $\bu(\bx,t)$ to denote the set of minimizers, which may be not a singleton set when the minimizer is not unique. We can write
\be\label{eqt: minplus_singletime_minimizer}
\begin{split}
\bu(\bx,t) &= \argmin_{\bu\in\R^n} \left\{ \min_{i \in \{1,\dots,m\}} J_i(\bu) + t H^*\left(\frac{\bx - \bu}{t}\right)\right\}\\
&= \argmin_{\bu\in\R^n} \min_{i \in \{1,\dots,m\}}\left\{  J_i(\bu) + t H^*\left(\frac{\bx - \bu}{t}\right)\right\}\\
&= \bigcup_{i\in I(\bx,t)}\argmin_{\bu\in\R^n} \left\{ J_i(\bu) + t H^*\left(\frac{\bx - \bu}{t}\right)\right\},
\end{split}
\ee
where the index set $I(\bx,t)$ is defined by 
\be\label{eqt: I_singletime}
I(\bx,t) = \argmin_{i \in \{1,\dots,m\}} S_i(\bx,t).
\ee

\def\edgeset {E}
\def\edgesubset {\Omega}
A specific example is when the regularization term $J$ is the truncated regularization term with pairwise interactions in the following form
\be \label{eqt: truncatedprior}
J(\bx) = \sum_{(i,j) \in \edgeset} w_{ij} f(x_i - x_j),  \text{ for each }\bx =(x_1, \dots, x_n)\in\R^n,
\ee
where $w_{ij} \geqslant 0$,  $f(x) = \min\{g(x), 1\}$ for some convex function $g\colon \R\to\R$ and $\edgeset=\{1,\dots,n\} \times \{1,\dots,n\}$. This function can be written as the minimum of a collection of convex functions $J_{\edgesubset}\colon \R^n\to \R$ as the following
\benn
J(\bx) = \min_{\edgesubset\subseteq \edgeset} J_{\edgesubset},
\eenn
with each $J_{\edgesubset}$ defined by
\benn
J_{\edgesubset} := \left\{\sum_{(i,j)\in \edgesubset} w_{ij} + \sum_{(i,j)\not\in \edgesubset} w_{ij}g(x_i- x_j) \right\},
\eenn
where $\Omega$ is any subset of $\edgeset$. The truncated regularization term \eqref{eqt: truncatedprior} can therefore be written in the form of \eqref{eqt: minJi}, and hence the minimizer to the corresponding optimization problem \eqref{eq:canonical-optimization-pb} with the non-convex regularization term $J$ in \eqref{eqt: truncatedprior} can be computed using \eqref{eqt: minplus_singletime_minimizer}.

We give here two examples of truncated regularization term with pairwise interactions in the form of~\eqref{eqt: truncatedprior}. First, let $g$ be the $\ell^1$ norm. Then $J$ is the truncated discrete Total Variation regularization term defined by
\be \label{eqt: truncatedTV}
J(\bx) = \sum_{(i,j)\in \edgeset} w_{ij}\min\{|x_i - x_j|,1\},  \text{ for each }\bx =(x_1, \dots, x_n)\in\R^n.
\ee
This function $J$ can be written as the formula \eqref{eqt: truncatedprior} with $f\colon \R\to \R$ given by Eq.~\eqref{eqt:l1}. Second, let $g$ be the quadratic function. Then $J$ is the half-quadratic regularization term defined by
\begin{equation}\label{eqt:half_quad_reg}
J(\bx) = \sum_{(i,j)\in \edgeset} w_{ij}\min\{(x_i - x_j)^2,1\},  \text{ for each }\bx =(x_1, \dots, x_n)\in\R^n.
\end{equation}
This function $J$ can be written as the formula \eqref{eqt: truncatedprior} with $f\colon \R\to \R$ given by Eq.~\eqref{eqt:half_quadratic}. This specific form of edge-preserving prior was investigated in the seminal works of~\cite{charbonnieretal.97.ip,geman_yang.95.ip,gemanReynolds.92.pami}. Several algorithms have been proposed to solve
the resultant non-convex optimization problem~\eqref{eqt: minplus_singletime}, i.e., the solution to the corresponding HJ PDE, for some specific choice of data fidelity terms (e.g.,~\cite{allain2006global,didier2001convex,geman_yang.95.ip,gemanReynolds.92.pami,nikolova2005analysis,champagnat2004connection,nikolova2001fast}).

Suppose now, for general regularization terms $J$ in the form of~\eqref{eqt: truncatedprior}, that we have Gaussian noise. Then the data fidelity term is quadratic and $H(\bp) = \frac{1}{2}\|\bp\|_2^2$ and $t = \frac{1}{\lambda}$. Hence, for this example, using \eqref{eqt: minplus_singletime_minimizer}, we obtain the set of minimizers
\benn
\begin{split}
\bu(\bx,t) &= \bigcup_{\edgesubset\in I(\bx,t)}\argmin_{\bu\in\R^n} \left\{ J_{\edgesubset}(\bu) + t H^*\left(\frac{\bx - \bu}{t}\right)\right\}\\
&= \bigcup_{\edgesubset\in I(\bx,t)}\argmin_{\bu\in\R^n} \left\{ \sum_{(i,j)\not\in \edgesubset} w_{ij} g(u_i- u_j) + \frac{1}{2t} \|\bx - \bu\|_2^2\right\}\\
&= \bigcup_{\edgesubset\in I(\bx,t)}\{\bx - t\nabla_{\bx} S_{\edgesubset}(\bx,t)\}
\end{split}
\eenn
where 
\benn
S_{\edgesubset}(\bx,t) = \sum_{(i,j)\in \edgesubset} w_{ij} + \min_{\bu\in\R^n} \left\{ \sum_{(i,j)\not\in \edgesubset} w_{ij}g(u_i- u_j) + \frac{1}{2t} \|\bx - \bu\|_2^2\right\}
\eenn
and 
\benn
I(\bx,t) = \argmin_{\edgesubset\subseteq \edgeset} S_{\edgesubset}(\bx,t).
\eenn

\bigbreak
The same result also holds for the multi-time HJ PDE system \eqref{eqt: multitime_HJ}. Indeed, if $J$ is a non-convex regularization term given by \eqref{eqt: minJi}, and $S, S_j\colon \R^n\times (0,+\infty)^N\to \R$ are the solutions to the multi-time HJ PDE system \eqref{eqt: multitime_HJ} with initial data $J$ and $J_i$, respectively, then similarly we have the min-plus linearity of the semi-group under certain assumptions. Specifically, if the Lax-Oleinik formula \eqref{eqt: multitime_LO} solves the multi-time HJ PDE system \eqref{eqt: multitime_HJ} for each $i\in \{1,\dots, m\}$ (for instance, when $H$ and $J_i$ satisfy the assumptions in Prop.~\ref{prop: firstorder_2} for each $i\in\{1,\dots, m\}$), then there holds
\be \label{eqt: minplus_multitime}
\begin{split}
S(\bx, t_1, \dots, t_N) &= \min_{\bu_1,\dots, \bu_N\in\R^n} \left\{ \min_{i \in \{1,\dots,m\}} J_i\left(\bx - \sum_{j=1}^N\bu_j\right) + \sum_{j=1}^N t_j H_j^*\left(\frac{\bu_j}{t_j}\right)\right\}\\
&=  \min_{i \in \{1,\dots,m\}}\left\{ \min_{\bu_1,\dots, \bu_N\in\R^n}\left\{ J_i\left(\bx - \sum_{j=1}^N\bu_j\right) + \sum_{j=1}^N t_j H_j^*\left(\frac{\bu_j}{t_j}\right)\right\}\right\}\\
&= \min_{i \in \{1,\dots,m\}} S_i(\bx,t_1, \dots, t_N).
\end{split}
\ee
Let $M\subset \R^{n\times N}$ be the set of minimizers of~\eqref{eqt: multitime_LO} with $J$ given by~\eqref{eqt: minJi}.
Then $M$ satisfies
\be\label{eqt: minplus_multitime_minimizer}
\begin{split}
M
&= \argmin_{\bu_1,\dots, \bu_N\in\R^n} \left\{ \min_{i \in \{1,\dots,m\}} J_i\left(\bx - \sum_{j=1}^N\bu_j\right) + \sum_{j=1}^N t_j H_j^*\left(\frac{\bu_j}{t_j}\right)\right\}\\
&= \bigcup_{i\in I(\bx,t_1,\dots, t_N)}\argmin_{\bu_1,\dots, \bu_N\in\R^n} \left\{ J_i\left(\bx - \sum_{j=1}^N\bu_j\right) + \sum_{j=1}^N t_j H_j^*\left(\frac{\bu_j}{t_j}\right)\right\},
\end{split}
\ee
where the index set $I(\bx,t_1,\dots, t_N)$ is defined by
\be\label{eqt: I_multitime}
I(\bx,t_1,\dots, t_N) = \argmin_{i \in \{1,\dots,m\}} S_i(\bx,t_1,\dots, t_N).
\ee
As a result, we can use \eqref{eqt: minplus_multitime_minimizer} to obtain the minimizers of the decomposition model \eqref{eqt: 1order_decomposition} with the non-convex regularization term $J$ in the form of \eqref{eqt: minJi}, such as the function in \eqref{eqt: truncatedprior} and the truncated Total Variation function \eqref{eqt: truncatedTV}.

\bigbreak
In summary, one can compute the minimizers of the optimization problems \eqref{eq:canonical-optimization-pb} and \eqref{eqt: 1order_decomposition} with a non-convex function $J$ in the form of \eqref{eqt: minJi} using the aforementioned min-plus algebra technique. Furthermore, this technique can be extended to handle other cases. For instance, in the denoising model \eqref{eq:canonical-optimization-pb}, if the data fidelity term $D$ is in the form of \eqref{eqt: minJi} and the prior term $\frac{J(\bu)}{\lambda}$ can be written as $tH^*\left(\frac{\bu}{t}\right)$, then one can still compute the minimizer of this problem using the min-plus algebra technique on the HJ PDE with initial data $D$. Similarly, because of the symmetry in the decomposition model \eqref{eqt: 1order_decomposition}, if there is only one non-convex term $f_j$ and if it can be written in the form of \eqref{eqt: minJi}, then one can apply the min-plus algebra technique to the multi-time HJ PDE with initial data $f_j$.

In general, however, there is a drawback to the min-plus algebra technique. To compute the minimizers using \eqref{eqt: minplus_singletime_minimizer} and \eqref{eqt: minplus_multitime_minimizer}, we need to compute the index set $I(\bx,t)$ and $I(\bx, t_1, \dots, t_N)$ defined in \eqref{eqt: I_singletime} and \eqref{eqt: I_multitime}, which involves solving $m$ HJ PDEs to obtain the solutions $S_1, \dots, S_m$. When $m$ is too large, this approach is impractical since it involves solving too many  HJ PDEs. 
For instance, if $J$ is the truncated Total Variation in \eqref{eqt: truncatedTV}, the number $m$ equals the number of subsets of the set $\edgeset$, i.e., $m=2^{|\edgeset|}$, which is computationally intractable. Hence, in general, it is impractical to use \eqref{eqt: minplus_singletime_minimizer} and \eqref{eqt: minplus_multitime_minimizer} to solve the problems \eqref{eq:canonical-optimization-pb} and \eqref{eqt: 1order_decomposition} where the regularization term $J$ is given by the truncated Total Variation. The same issue arises when the truncated Total Variation is replaced by half-quadratic regularization. Several authors attempted to address this intractability for half-quadratic regularizations by proposing heuristic optimization methods that aim to compute a global minimizer~\cite{allain2006global,didier2001convex,geman_yang.95.ip,gemanReynolds.92.pami,nikolova2005analysis,champagnat2004connection,nikolova2001fast}.

\begin{figure}[htbp]
\centering
\includegraphics[width = 0.5\textwidth]{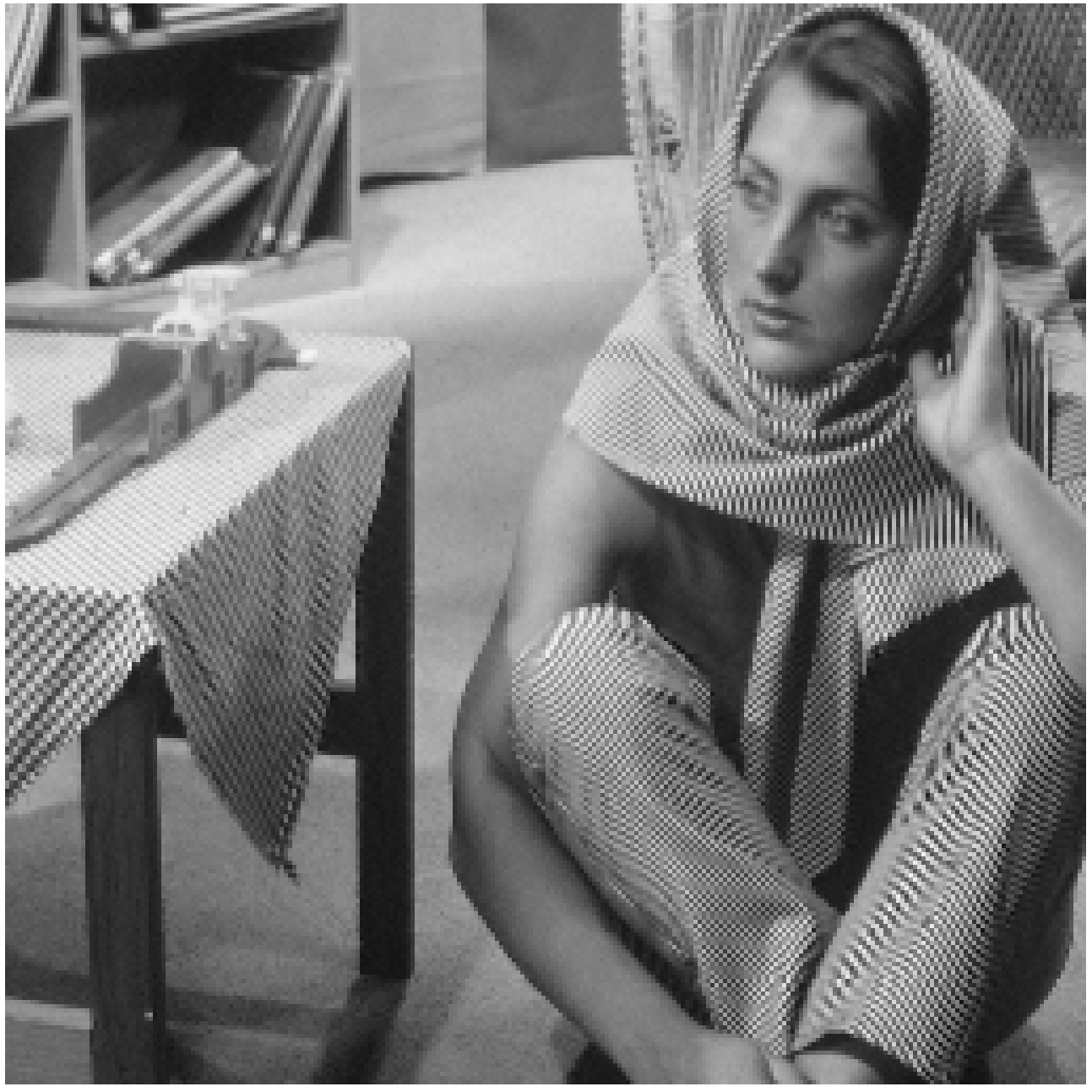}
\caption{The input image $\bx$ (``Barbara'') in the example in Sec.~\ref{sec: eg}.}
\label{fig:barbara}
\end{figure}

\begin{figure}[htbp]
\centering
\subfloat[]{\label{fig:y1}\includegraphics[width = 0.5\textwidth]{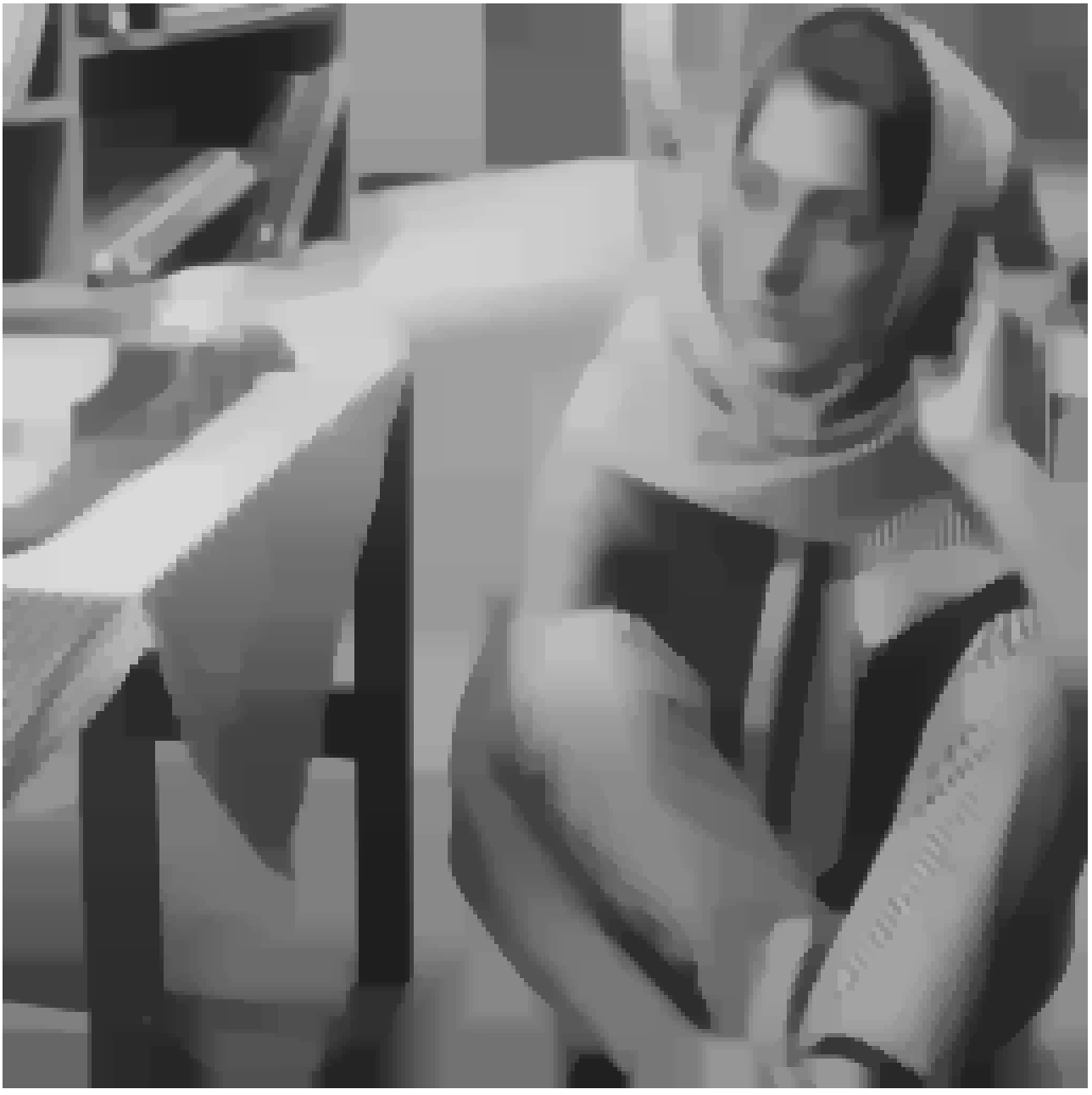}}\\
\subfloat[]{\label{fig:z1}\includegraphics[width = 0.5\textwidth]{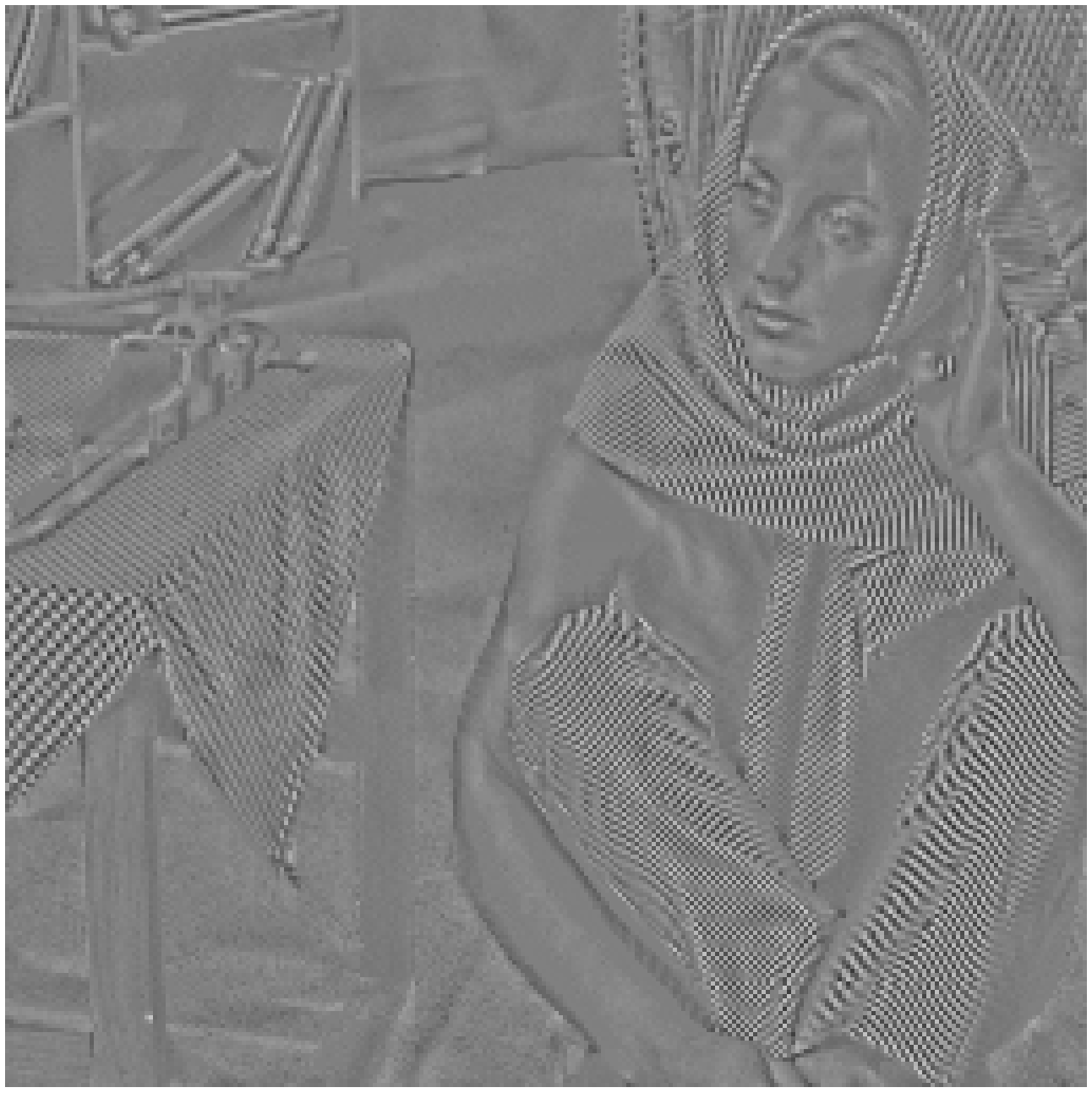}}
\subfloat[]{\label{fig:w1}\includegraphics[width = 0.5\textwidth]{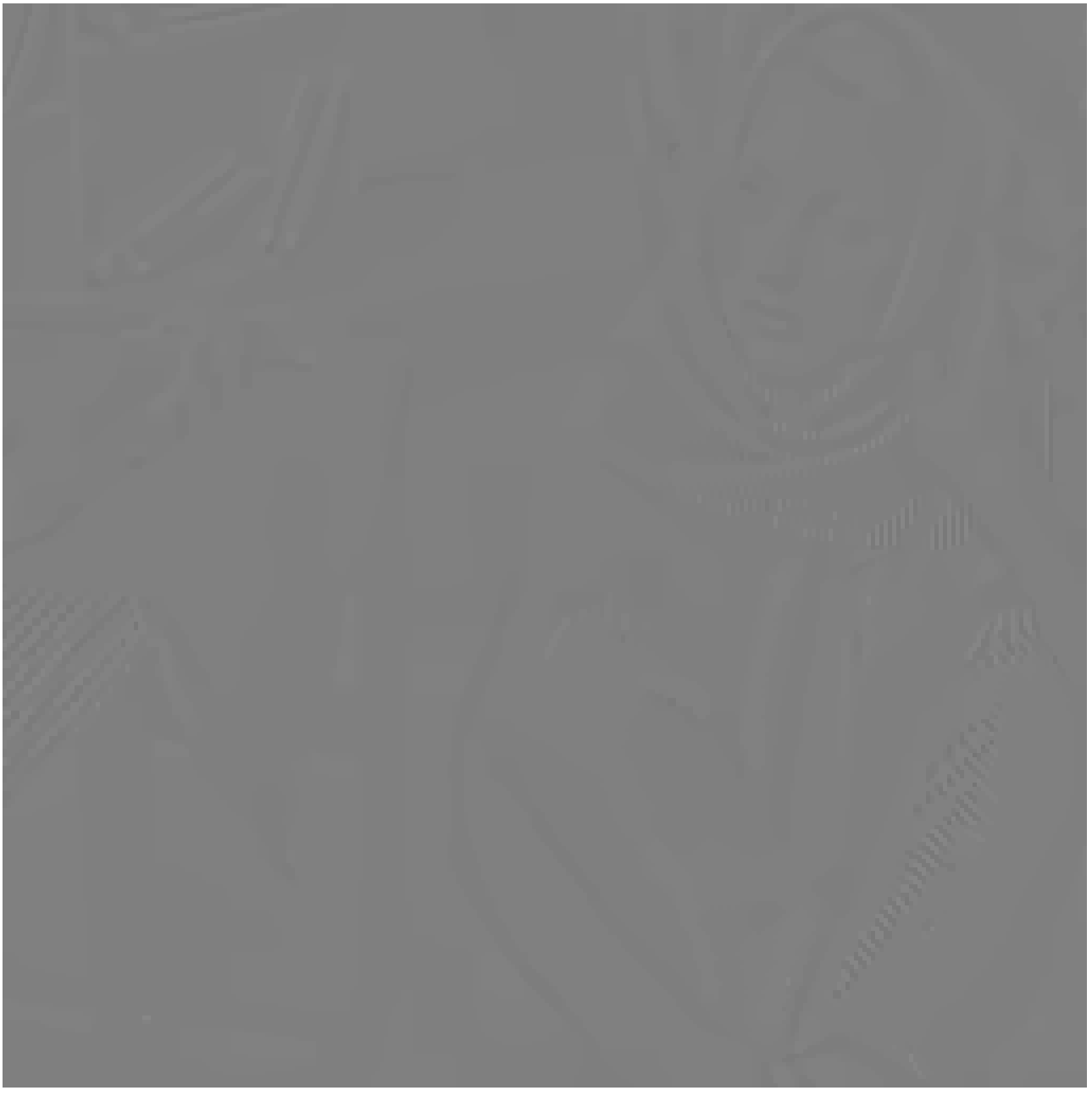}}
\caption{The minimizer of the first problem in~\eqref{eqt: eg_defS1S2}. The output images $\bx-\bu_1-\bu_2$, $\bu_1 + 0.5$ and $\bu_2+0.5$ are shown in (a), (b) and (c), respectively.}
\label{fig:1}
\end{figure}

\begin{figure}[htbp]
\centering
\subfloat[]{\label{fig:y2}\includegraphics[width = 0.5\textwidth]{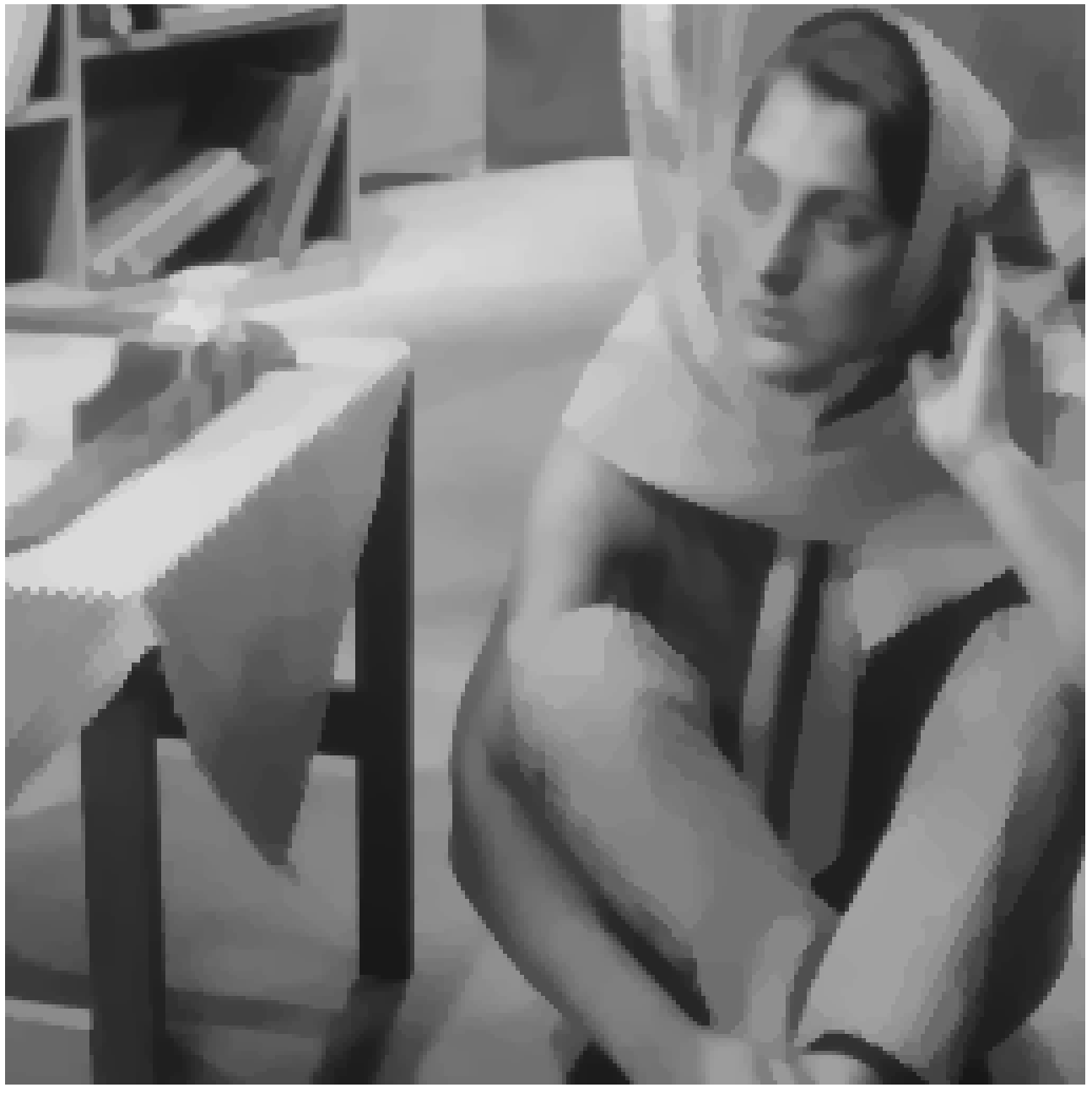}}\\
\subfloat[]{\label{fig:z2}\includegraphics[width = 0.5\textwidth]{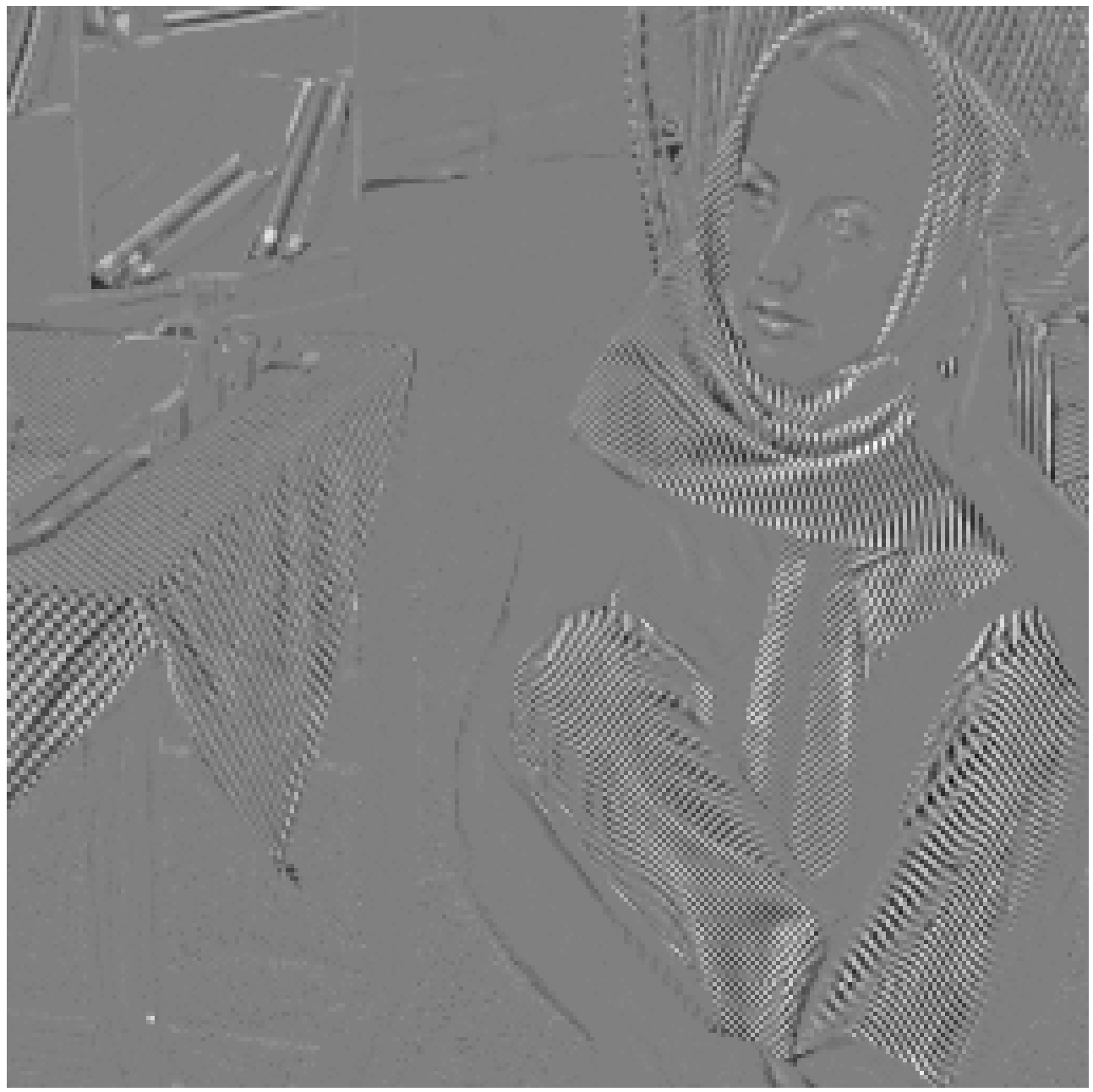}}
\subfloat[]{\label{fig:w2}\includegraphics[width = 0.5\textwidth]{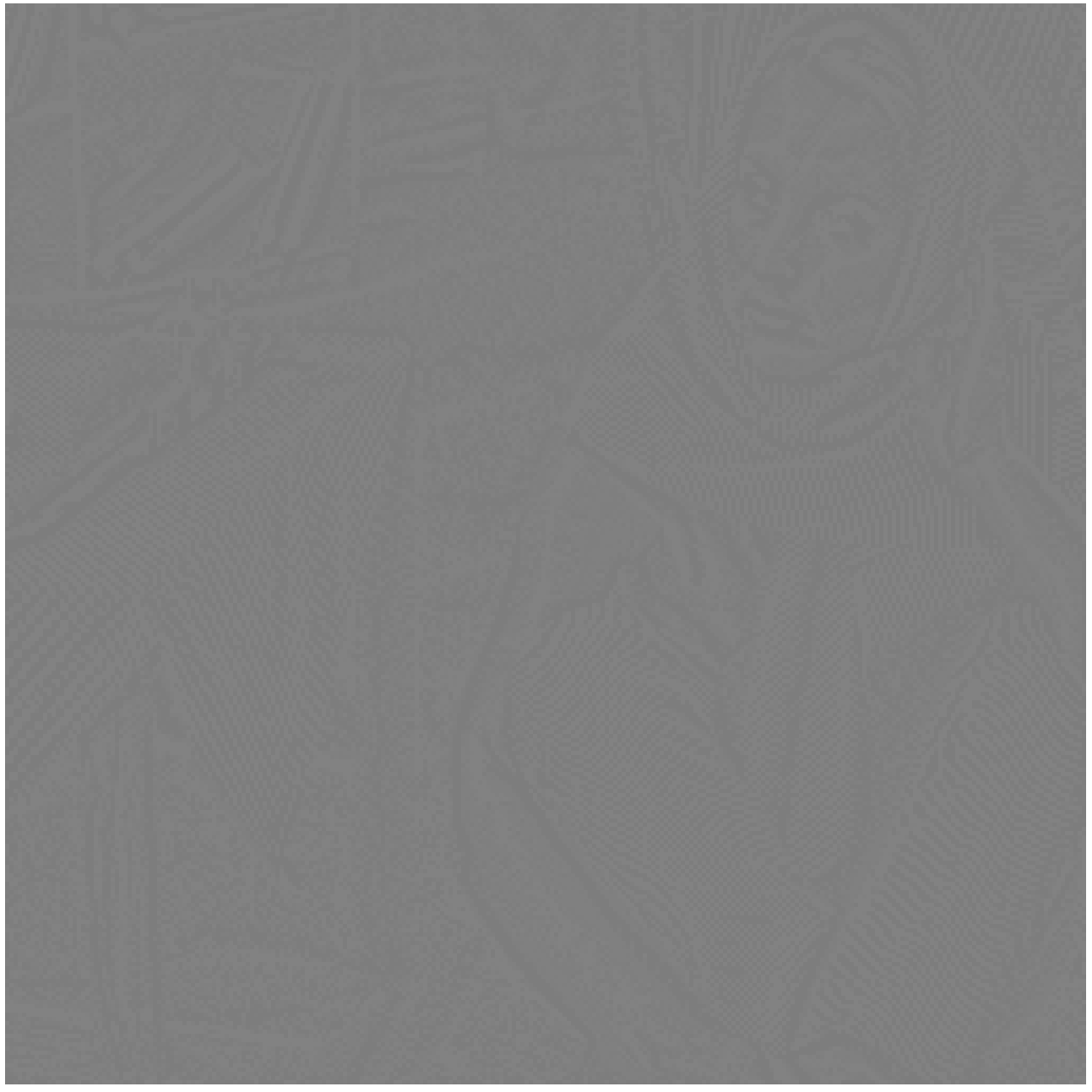}}
\caption{The minimizer of the second problem in~\eqref{eqt: eg_defS1S2}. The output images $\bx-\bv_1-\bv_2$, $\bv_1 + 0.5$ and $\bv_2+0.5$ are shown in (a), (b) and (c), respectively.}
\label{fig:2}
\end{figure}

\subsection{Application to certain decomposition problems}\label{sec: eg}
In this section, we demonstrate how to use our formulation described in the previous sections to solve certain image decomposition problems. The variational formulation for image decomposition problems is in the form of~\eqref{eqt: 1order_decomposition}, where the input image $\bx\in\Rn$ is decomposed into three components, which includes the geometrical part $\bx-\bu_1-\bu_2$, the texture part $\bu_1$, and the noise $\bu_2$. 
The regularization function $J$ for the geometrical part $\bx-\bu_1-\bu_2$ is chosen to be the widely used Total Variation regularization function in order to preserve edges in the image. Here, we use the anisotropic Total Variation semi-norm (see e.g, \cite{darbon2006image,darbon2015convex}) denoted by $|\cdot|_{TV}$. The noise is assumed to be Gaussian, and hence the data fidelity term $f_2$ is set to be the quadratic function. Many texture models have been proposed (see~\cite{aujol2003image, aujol2005image,LetGuen2014Cartoon,winkler.03.book} and the references in these papers).
For instance, the indicator function of the unit ball with respect to Meyer's norm is used in~\cite{aujol2003image,aujol2005image}, and the $\ell^1$ norm is used in~\cite{LetGuen2014Cartoon}. Note that each texture model has some pros and cons and, to our knowledge, it remains an open problem whether one specific texture model is better than the others.  
In this example, we combine different texture regularizations proposed in the literature by
taking the minimum of the indicator function of the unit ball with respect to Meyer's norm and the $\ell^1$ norm. In other words, we consider the following variational problem
\begin{equation}\label{eqt: eg_decomp}
\min_{\bu_1,\bu_2\in\R^n} \left\{ J(\bx - \bu_1 - \bu_2) + t_1 g\left(\frac{\bu_1}{t_1}\right) +  \frac{1}{2t_2}\|\bu_2\|_2^2\right\},
\end{equation}
where $J\colon \Rn\to\R$ and $g\colon\Rn\to\R\cup\{+\infty\}$ are defined by
\begin{equation*}
    J(\by) := |\by|_{TV}, \quad\quad
    g(\by):= \min\{J^*(\by), \|\by\|_1\},
\end{equation*}
for each $\by\in\Rn$. Problem~\eqref{eqt: eg_decomp} is equivalent to the following mixed discrete-continuous optimization problem
\begin{equation}\label{eqt: eg_dis_cont}
    \min_{\bu_1,\bu_2\in\R^n}\min_{k\in\{1,2\}} \left\{ J(\bx - \bu_1 - \bu_2) + t_1 g_k\left(\frac{\bu_1}{t_1}\right) +  \frac{1}{2t_2}\|\bu_2\|_2^2\right\},
\end{equation}
where $g_1(\by):= J^*(\by)$ and $g_2(\by) := \|\by\|_1$ for each $\by\in\Rn$.
Note that solving mixed discrete-continuous optimization is hard in general (see~\cite{encyclopediaOptimization09} for instance). However, we shall see that our proposed approach yields efficient optimization algorithms.
Since the function $g$ is the minimum of two convex functions, the problem~\eqref{eqt: eg_decomp}
fits into our formulation, and can be solved using a similar idea as in~\eqref{eqt: minplus_multitime} and~\eqref{eqt: minplus_multitime_minimizer}. 
To be specific, define the two functions $S_1$ and $S_2$ by
\begin{equation}\label{eqt: eg_defS1S2}
    \begin{split}
    S_1(\bx,t_1,t_2)&:= \min_{\bu_1,\bu_2\in\R^n} \left\{ J(\bx - \bu_1 - \bu_2) + t_1J^*\left(\frac{\bu_1}{t_1}\right) +  \frac{1}{2t_2}\|\bu_2\|_2^2\right\},\\
    S_2(\bx,t_1,t_2)&:= \min_{\bu_1,\bu_2\in\R^n} \left\{ J(\bx - \bu_1 - \bu_2) + \|\bu_1\|_1 +  \frac{1}{2t_2}\|\bu_2\|_2^2\right\},
    \end{split}
\end{equation}
where the sets of the minimizers in the two minimization problems above are denoted by $M_1(\bx,t_1,t_2)$ and $M_2(\bx,t_1,t_2)$, respectively.
Using a similar argument as in~\eqref{eqt: minplus_multitime} and~\eqref{eqt: minplus_multitime_minimizer}, we conclude that the minimal value in~\eqref{eqt: eg_decomp} equals $\min\{S_1(\bx,t_1,t_2), S_2(\bx,t_1,t_2)\}$, and the set of minimizers in~\eqref{eqt: eg_decomp}, denoted by $M(\bx,t_1,t_2)$, satisfies
\begin{equation}\label{eqt: eg_minimizer}
    M(\bx,t_1,t_2) = \begin{dcases}
    M_1(\bx,t_1,t_2) & S_1(\bx,t_1,t_2) < S_2(\bx,t_1,t_2),\\
    M_2(\bx,t_1,t_2) & S_1(\bx,t_1,t_2) > S_2(\bx,t_1,t_2),\\
    M_1(\bx,t_1,t_2)\cup M_2(\bx,t_1,t_2)  & S_1(\bx,t_1,t_2) = S_2(\bx,t_1,t_2).
    \end{dcases}
\end{equation}
As a result, we solve the two minimization problems in~\eqref{eqt: eg_defS1S2} first, and then obtain the minimizers using~\eqref{eqt: eg_minimizer} by comparing the minimal values $S_1(\bx,t_1,t_2)$ and $S_2(\bx,t_1,t_2)$.

Here, we present a numerical result. 
We solve the first optimization problem in~\eqref{eqt: eg_defS1S2} by a splitting method, where each subproblem can be solved using the proximal operator of the anisotropic Total Variation (for more details, see \cite{darbon2019decomposition}). Similarly, a splitting method is used to split the second optimization problem in~\eqref{eqt: eg_defS1S2} to two subproblems, which are solved using the proximal operators of the anisotropic Total Variation and the $\ell^1$-norm, respectively. 
To compute the proximal point of the anisotropic Total Variation, the algorithm in~\cite{chambolle.09.ijcv,darbon2006image,hochbaum.01.jacm} is adopted and it computes the proximal point without numerical errors.
The input image $\bx$ is the image ``Barbara'' shown in Fig.~\ref{fig:barbara}. The parameters are set to be $t_1 = 0.07$ and $t_2 = 0.01$. Let $(\bu_1,\bu_2)\in M_1(\bx,t_1,t_2)$ and $(\bv_1,\bv_2)\in M_2(\bx,t_1,t_2)$ be respectively the minimizers of the two minimization problems in~\eqref{eqt: eg_defS1S2} solved by the aforementioned splitting methods. We show these minimizers and the related images in Figs.~\ref{fig:1} and~\ref{fig:2}.
To be specific, the decomposition components $\bx-\bu_1-\bu_2$, $\bu_1+0.5$ and $\bu_2+0.5$ given by the first optimization problem in~\eqref{eqt: eg_defS1S2} are shown in Figs.~\ref{fig:y1},~\ref{fig:z1} and~\ref{fig:w1}, respectively. The decomposition components $\bx-\bv_1-\bv_2$, $\bv_1 +0.5$ and $\bv_2+0.5$ given by the second optimization problem in~\eqref{eqt: eg_defS1S2} are shown in Figs.~\ref{fig:y2},~\ref{fig:z2} and~\ref{fig:w2}, respectively. We also compute the optimal values $S_1(\bx,t_1,t_2)$ and $S_2(\bx,t_1,t_2)$, and obtain
\begin{equation*}
    S_1(\bx,t_1,t_2) = 1832.81, \quad\quad
    S_2(\bx,t_1,t_2) = 4171.33.
\end{equation*}
Since $S_1(\bx,t_1,t_2) < S_2(\bx,t_1,t_2)$, we conclude that $(\bu_1, \bu_2)$ is a minimizer in the decomposition problem~\eqref{eqt: eg_decomp}, and the minimal value equals $1832.81$. In other words, the optimal decomposition given by~\eqref{eqt: eg_decomp} is shown in Fig.~\ref{fig:1}.

%% file: secondorder.tex
\section{Viscous Hamilton--Jacobi PDEs and Bayesian estimation}
\label{sec:viscousHJ}

In contrast to variational approaches that frame imaging problems as optimization problems, Bayesian approaches frame them in a probabilistic framework. This framework combines observed data through a likelihood function (which models the noise corrupting the unknown image) and prior knowledge through a prior distribution (which models known properties of the image to reconstruct) to generate a posterior distribution from which an appropriate decision rule can select a meaningful image estimate. In this section, we present an analogue of the min-plus algebra technique discussed in Sect.~\ref{subsec:minplusMultiTimeHJ} for certain Bayesian posterior mean estimators.

\subsection{Viscous HJ PDEs and posterior mean estimators for log-concave models}
\label{subsec:visousHJlogconcave}

Consider the following class of Bayesian posterior distributions
\begin{equation}\label{eqt: post_dist}
q(\bu|(\bx,t,\epsilon)) \coloneqq \frac{e^{-\left(J(\bu) + \frac{1}{2t}\normsq{\bx-\bu}\right)/\epsilon}}{\int_{\Rn} e^{-\left(J(\bu) + \frac{1}{2t}\normsq{\bx-\bu}\right)/\epsilon} \diff\bu},
\end{equation}
where $\bx \in \Rn$ is the observed image with $n$ pixels, and $t$ and $\epsilon$ are positive parameters. The posterior distribution~\eqref{eqt: post_dist} is proportional to the product of a log-concave prior $\bu \mapsto e^{-J(\bu)/\epsilon}$ (possibly improper) and a Gaussian likelihood function $\bu \mapsto e^{-\frac{1}{2t\epsilon}\normsq{\bx-\bu}}$. This class of posterior distributions generates the family of Bayesian posterior mean estimators $\bu_{PM}\colon \Rn\times(0,+\infty)\times(0,+\infty) \to \Rn$ defined by
\begin{equation}\label{eqt: pm_estimator}
    \bu_{PM}(\bx,t,\epsilon) \coloneqq \int_{\Rn} \bu\,q(\bu|(\bx,t,\epsilon)) \diff\bu.
\end{equation}
These are Bayesian estimators because they minimize the mean squared error (\cite{kay1993fundamentals}, pages 344-345):
\begin{equation}\label{eqt: MSE_def}
\bu_{PM}(\bx,t,\epsilon) = \argmin_{\bu \in \Rn} \int_{\Rn} \normsq{\bar{\bu}-\bu}\, q(\bar{\bu}|(\bx,t,\epsilon)) \diff\bar{\bu}.
\end{equation}
They are frequently called minimum mean squared error estimators for this reason.

The class of posterior distributions~\eqref{eqt: post_dist} also generates the family of maximum a posteriori estimators $\bu_{MAP}:\Rn\times(0,+\infty) \to \Rn$ defined by
\begin{equation}\label{eqt: map_estimator}
    \bu_{MAP}(\bx,t) = \argmin_{\bu \in \Rn}\left\{J(\bu) + \frac{1}{2t}\normsq{\bx-\bu}\right\},
\end{equation}
where $\bu_{MAP}(\bx,t)$ is the mode of the posterior distribution~\eqref{eqt: post_dist}. Note that the MAP estimator is also the minimizer of the solution~\eqref{eqt: 1order_LO} to the first order HJ PDE~\eqref{eqt: 1order_HJ} with Hamiltonian $H = \frac{1}{2}\normsq{\cdot}$ and initial data $J$.

There is a large body of literature on posterior mean estimators for image restoration problems~(see e.g., \cite{demoment1989image,kay1993fundamentals,winkler.03.book}). In particular, original connections between variational problems and Bayesian methods have been investigated in~\cite{louchet2008modeles,louchet2013posterior,burger2014maximum,burger2016bregman,gribonval2011should,gribonval2013reconciling,gribonval2018bayesian,darbon2020bayesian}. In particular, in \cite{darbon2020bayesian}, the authors described original connections between Bayesian posterior mean estimators and viscous HJ PDEs when $J \in \gmRn$ and the data fidelity term is Gaussian. We now briefly describe these connections here. 

Consider the function $S_\epsilon: \Rn \times (0,+\infty) \to \R$ defined by
\begin{equation}\label{eqt: Seps_formula}
    S_\epsilon(\bx,t) = -\epsilon\ln\left(\frac{1}{(2\pi t\epsilon)^{n/2}} \int_{\Rn} e^{-\left(J(\bu) + \frac{1}{2t}\normsq{\bx-\bu}\right)/\epsilon} \diff\bu\right),
\end{equation}
which is proportional to the negative logarithm of the partition function of the posterior distribution~\eqref{eqt: post_dist}. Under appropriate assumptions on the regularization term $J$ (see proposition~\ref{prop: secondorder_1} below), formula~\eqref{eqt: Seps_formula} corresponds to a Cole-Hopf transform \cite{evans1998partial} and is the solution to the following viscous HJ PDE
\begin{equation}\label{eqt: 2ndorder_HJ}
\begin{dcases}
\frac{\partial S_\epsilon}{\partial t}(\bx,t) + \frac{1}{2}\normsq{\nabla_{\bx}S_\epsilon(\bx,t)} = \frac{\epsilon}{2}\Laplacian S_\epsilon(\bx,t) &\bx\in \Rn, t>0,\\
S_\epsilon(\bx,0) = J(\bx) &\bx\in \Rn,
\end{dcases}
\end{equation}
where $J$ is the initial data. The solution to this PDE is also related to the first-order HJ PDE~\eqref{eqt: 1order_HJ} when the Hamiltonian is $H = \frac{1}{2}\normsq{\cdot}$. The following proposition, which is given in~\cite{darbon2020bayesian}, describes these connections.
\begin{proposition}\label{prop: secondorder_1}
Assume $J \in \gmRn$, $\inter{(\dom J)} \neq \varnothing$, and $\inf_{\bu \in \Rn} J(\bu) = 0$. Then for every $\epsilon > 0$, the unique smooth solution $S_\epsilon\colon \Rn\times(0,+\infty)\to(0,+\infty)$ to the HJ PDE~\eqref{eqt: 2ndorder_HJ} is given by formula~\eqref{eqt: Seps_formula}, where $(\bx,t)\mapsto S_\epsilon(\bx,t) -\frac{n\epsilon}{2}\ln t$ is jointly convex. Moreover, for each $\bx \in \Rn$, $t>0$, and $\epsilon > 0$, the posterior mean estimator~\eqref{eqt: pm_estimator} and minimum mean squared error in~\eqref{eqt: MSE_def} (with $\bu = \bu_{PM}(\bx,t,\epsilon)$) satisfy, respectively, the formulas
\begin{equation}\label{eqt: prop1_relation_bayes_1}
\bu_{PM}(\bx,t,\epsilon) = \bx - t\nabla_{\bx}S_\epsilon(\bx,t)
\end{equation}
and
\begin{equation} \label{eqt: prop1_relation_bayes_2}
    \int_{\Rn} \normsq{\bu_{PM}(\bx,t,\epsilon)-\bu} \, q(\bu|(\bx,t,\epsilon)) \diff\bu = nt\epsilon - t^2\epsilon\Laplacian S_{\epsilon}(\bx,t).
\end{equation}
In addition, for every $\bx \in \Rn$ and $t>0$, the limits of $\lim_{\substack{\epsilon\to0\\\epsilon>0}} S_\epsilon(\bx,t)$ and $\lim_{\substack{\epsilon\to0\\\epsilon>0}} \bu_{PM}(\bx,t,\epsilon)$ exist and converge uniformly over every compact set of $\Rn\times(0,+\infty)$ in $(\bx,t)$. Specifically, we have
\begin{equation} \label{eqt: eps_lim_seps}
    \lim_{\substack{\epsilon\to0\\\epsilon>0}} S_\epsilon(\bx,t) = \min_{\bu\in\R^n} \left\{ J(\bu) + \frac{1}{2t}\normsq{\bx-\bu}\right\},
\end{equation}
where the right hand side solves uniquely the first order HJ PDE~\eqref{eqt: 1order_HJ} with Hamiltonian $H = \frac{1}{2}\normsq{\cdot}$ and initial data $J$, and
\begin{equation} \label{eqt: eps_lim_upm}
    \lim_{\substack{\epsilon\to0\\\epsilon>0}} \bu_{PM}(\bx,t,\epsilon) = \argmin_{\bu\in\R^n} \left\{ J(\bu) + \frac{1}{2t}\normsq{\bx-\bu}\right\}.
\end{equation}
\end{proposition}

Under convexity assumptions on $J$, the representation formulas \eqref{eqt: prop1_relation_bayes_1} and \eqref{eqt: prop1_relation_bayes_2} relate the posterior mean estimate and the minimum mean squared error to the spatial gradient and Laplacian of the solution to the viscous HJ PDE~\eqref{eqt: 2ndorder_HJ}, respectively. Hence one can compute the posterior mean estimator and minimum mean square error using the spatial gradient $\nabla_{\bx}S_\epsilon(\bx,t)$ and the Laplacian $\Laplacian S_{\epsilon}(\bx,t)$ of the solution to the HJ PDE~\eqref{eqt: 2ndorder_HJ}, respectively, or vice versa by computing the posterior mean and minimum mean square error using, for instance, Markov Chain Monte Carlo sampling strategies.

The limit~\eqref{eqt: eps_lim_upm} shows that the posterior mean $\bu_{PM}(\bx,t,\epsilon)$ converges to the maximum a posteriori $\bu_{MAP}(\bx,t)$ as the parameter $\epsilon \to 0$. A rough estimate of the squared Euclidean distance between the posterior mean estimator~\eqref{eqt: pm_estimator} and the maximum a posteriori~\eqref{eqt: map_estimator} in terms of the parameters $t$ and $\epsilon$ is given by
\begin{equation}\label{eqt: diff_pm_min}
    \normsq{\bu_{PM}(\bx,t,\epsilon) - \bu_{MAP}(\bx,t)} \leqslant nt\epsilon.
\end{equation}

\subsection{On viscous HJ PDEs with certain non log-concave priors}

So far, we have assumed that the regularization term $J$ in the posterior distribution~\eqref{eqt: post_dist} and Proposition~\ref{prop: secondorder_1} is convex. Here, we consider an analogue of the min-plus algebra technique designed for certain first order HJ PDEs tailed to viscous HJ PDEs, which will enable us to derive representation formulas for posterior mean estimators of the form of~\eqref{eqt: pm_estimator} whose priors are sums of log-concave priors, i.e., to certain mixture distributions.

Remember that the min-plus algebra technique for first order HJ PDEs described in Sect.~\ref{subsec:minplusMultiTimeHJ} involves initial data of the form $\min_{i\in\{1,\dots,m\}} J_i(\bx)$ where each $J_i\colon\Rn\to\R\cup\{+\infty\}$ is convex. Consider now initial data of the form
\begin{equation}\label{eqt: smoothminJi}
\displaystyle
    J(\bx) = -\epsilon\ln\left(\sum_{i=1}^{m}e^{-J_i(\bx)/\epsilon}\right).
\end{equation}
Note that formula~\eqref{eqt: smoothminJi} approximates the non-convex term~\eqref{eqt: minJi} in that
\[
\lim_{\substack{\epsilon\to0\\\epsilon>0}} -\epsilon\ln\left(\sum_{i=1}^{m}e^{-J_i(\bx)/\epsilon}\right) = \min_{i\in\{1,\dots,m\}} J_i(\bx) \text{ for each } \bx\in\R^n. 
\]
Now, assume $\inter{(\dom J_i)} \neq \varnothing$ for each $i \in \{1,\dots,m\}$, and let
\[
S_{i,\epsilon}(\bx,t) = -\epsilon\ln\left(\frac{1}{(2\pi t\epsilon)^{n/2}} \int_{\Rn} e^{-\left(J_i(\bu) + \frac{1}{2t}\normsq{\bx-\bu}\right)/\epsilon} \diff\bu\right),
\]
and
\[
\bu_{i,PM}(\bx,t,\epsilon) = \frac{\int_{\Rn} \bu\, e^{-\left(J_i(\bu) + \frac{1}{2t}\normsq{\bx-\bu}\right)/\epsilon} \diff\bu}{\int_{\Rn} e^{-\left(J_i(\bu) + \frac{1}{2t}\normsq{\bx-\bu}\right)/\epsilon} \diff\bu}
\]
denote, respectively, the solution to the viscous HJ PDE~\eqref{eqt: 2ndorder_HJ} with initial data $J_i$ and its associated posterior mean. Then, a short calculation shows that for every $\epsilon > 0$, the function $S_\epsilon(\bx,t)\colon\Rn\times(0,+\infty)\to\R$ defined by
\begin{equation}\label{eqt: seps_multipleJi}
\begin{split}
    S_\epsilon(\bx,t) &= -\epsilon\ln\left(\sum_{i=1}^{m} \frac{1}{(2\pi t\epsilon)^{n/2}} \int_{\Rn} e^{-\left(J_i(\bu) + \frac{1}{2t}\normsq{\bx-\bu}\right)/\epsilon} \diff\bu\right)\\
    &= -\epsilon\ln\left(\sum_{i=1}^{m} e^{-S_{i,\epsilon}(\bx,t)/\epsilon}\right)
\end{split}
\end{equation}
is the unique smooth solution to the viscous HJ PDE~\eqref{eqt: 2ndorder_HJ} with initial data~\eqref{eqt: smoothminJi}. As stated in Sect.~\ref{subsec:visousHJlogconcave}, the posterior mean estimate $\bu_{PM}(\bx,t,\epsilon)$ is given by the representation formula
\begin{equation} \label{eqt: prop1_relation_bayes_1_mult}
\bu_{PM}(\bx,t,\epsilon) = \bx - t\nabla_{\bx}S_\epsilon(\bx,t),
\end{equation}
which can be expressed in terms of the solutions $S_{i,\epsilon}(\bx,t)$, their spatial gradients $\nabla_{\bx}S_{i,\epsilon}(\bx,t)$, and posterior mean estimates $\bu_{i,PM}(\bx,t,\epsilon)$ as the weighted sums
\begin{equation}\label{eqt: i_pm_estimator}
\begin{split}
    \bu_{PM}(\bx,t,\epsilon) &= \bx - t\left(\frac{\sum_{i=1}^{m} \nabla_{\bx}S_{i,\epsilon}(\bx,t)e^{-S_{i,\epsilon}(\bx,t)/\epsilon}}{\sum_{i=1}^{m}e^{-S_{i,\epsilon}(\bx,t)/\epsilon}}\right)\\
    &= \frac{\sum_{i=1}^{m} \bu_{i,PM}(\bx,t,\epsilon)e^{-S_{i,\epsilon}(\bx,t)/\epsilon}}{\sum_{i=1}^{m}e^{-S_{i,\epsilon}(\bx,t)/\epsilon}}.
\end{split}
\end{equation}

As an application of this result, we consider the problem of classifying a noisy image $\bx \in \Rn$ using a Gaussian mixture model \cite{duda2012pattern}: Suppose $J_i(\bu) = \frac{1}{2\sigma_{i}^{2}}\normsq{\bu - \boldsymbol{\mu}_i}$, where $\boldsymbol{\mu}_i \in \Rn$ and $\sigma_i > 0$. The regularized minimization problem~\eqref{eqt: minplus_singletime} with quadratic data fidelity term $H = \frac{1}{2}\normsq{\cdot}$ is given by
\begin{equation}
\begin{split}
S_0(\bx, t) &= \min_{\bu\in\R^n} \left\{\min_{i \in \{1,\dots,m\}}\left\{\frac{1}{2\sigma_{i}^{2}}\normsq{\bu - \boldsymbol{\mu}_i} + \frac{1}{2t}\normsq{\bx-\bu}\right\}\right\}\\
&=  \min_{i \in \{1,\dots,m\}}\left\{ \min_{\bu\in\R^n}\left\{ \frac{1}{2\sigma_{i}^{2}}\normsq{\bu - \boldsymbol{\mu}_i} + \frac{1}{2t}\normsq{\bx-\bu}\right\}\right\} \\
&= \min_{i \in \{1,\dots,m\}} \left\{ \frac{1}{2(\sigma_{i}^{2}+t)}\normsq{\bx - \boldsymbol{\mu}_i}\right\}.
\end{split}
\end{equation}
Letting $I(\bx,t) = \argmin_{i\in\{1,\dots,m\}} \left\{ \frac{1}{2(\sigma_{i}^{2}+t)}\normsq{\bx - \boldsymbol{\mu}_i}\right\}$, the MAP estimator is then the collection
\[
\bu_{MAP}(\bx,t) = \bigcup_{i\in I(\bx,t)}\left\{\frac{\sigma_{i}^{2}\bx +t\boldsymbol{\mu}_i}{\sigma_{i}^2+t}\right\}.
\]
Consider now the initial data~\eqref{eqt: smoothminJi}:
\[
J(\bu) = -\epsilon\ln\left(\sum_{i=1}^{m} e^{-\frac{1}{2\sigma_{i}^{2}\epsilon}\normsq{\bu - \boldsymbol{\mu}_i}}\right).
\]
The solution $S_\epsilon(\bx,t)$ to the viscous HJ PDE~\eqref{eqt: 2ndorder_HJ} with initial data $J(\bx)$ is given by formula~\eqref{eqt: seps_multipleJi}, which in this case can be computed analytically:
\begin{equation}
S_\epsilon(\bx,t) = -\epsilon\ln\left(\sum_{i=1}^{m}\left(\frac{\sigma_{i}^{2}}{\sigma_{i}^{2}  + t}\right)^{n/2} e^{-\frac{1}{2(\sigma_{i}^{2} + t)\epsilon}\normsq{\bx - \boldsymbol{\mu}_i}}\right).
\end{equation}
Since $e^{-S_{i,\epsilon}(\bx,t)/\epsilon} = \left(\frac{\sigma_{i}^{2}}{\sigma_{i}^{2}  + t}\right)^{n/2} e^{-\frac{1}{2(\sigma_{i}^{2} + t)\epsilon}\normsq{\bx - \boldsymbol{\mu}_i}}$, we can write the corresponding posterior mean estimator~\eqref{eqt: i_pm_estimator} using the representation formulas~\eqref{eqt: prop1_relation_bayes_1_mult} and~\eqref{eqt: i_pm_estimator}:
\begin{equation}
    \begin{split}
        \bu_{PM}(\bx,t,\epsilon) &= \bx - t\nabla_{\bx}S_\epsilon(\bx,t)\\
        &= \frac{\sum_{i=1}^{m} \left(\frac{\sigma_{i}^{2}\bx + t\boldsymbol{\mu}_i}{\sigma_{i}^{2} + t}\right)  \left(\frac{\sigma_{i}^{2}}{\sigma_{i}^{2}  + t}\right)^{n/2} e^{-\frac{1}{2(\sigma_{i}^{2} + t)\epsilon}\normsq{\bx - \boldsymbol{\mu}_i}}}{\sum_{i=1}^{m}  \left(\frac{\sigma_{i}^{2}}{\sigma_{i}^{2}  + t}\right)^{n/2} e^{-\frac{1}{2(\sigma_{i}^{2} + t)\epsilon}\normsq{\bx - \boldsymbol{\mu}_i}}}.
    \end{split}
\end{equation}

%% file: conclusion.tex
\section{Conclusion}
\label{sec:conclusion}

In this chapter, we reviewed the connections of single time HJ PDEs with image denoising models and the connections of multi-time HJ PDEs with image decomposition models under convexity assumptions. Specifically, under some assumptions, the minimizers of these optimization problems can be computed using the spatial gradient of the solution to the corresponding HJ PDEs. We also proposed a min-plus algebra technique to cope with certain non-convex regularization terms in imaging sciences problems. This suggests that certain non-convex optimization problem can be solved by computing several convex sub-problems. For instance, if the denoising model~\eqref{eq:canonical-optimization-pb} or the image decomposition model~\eqref{eqt: 1order_decomposition} involves a non-convex regularization term $J$ that can be expressed as the minimum of $m$ convex sub-problems in the form of~\eqref{eqt: minJi}, then the minimizer of these non-convex problems can be solved using formulas~\eqref{eqt: minplus_singletime_minimizer} and~\eqref{eqt: minplus_multitime_minimizer}. However, when $m$ in \eqref{eqt: minJi} is too large, it is generally impractical to solve~\eqref{eqt: minplus_singletime_minimizer} and~\eqref{eqt: minplus_multitime_minimizer} using this min-plus technique because it involves solving too many HJ PDEs. However, our formulation yields practical numerical optimization algorithms for certain image denoising and decomposition problems.

We also reviewed connections between viscous HJ PDEs and a class of Bayesian methods and posterior mean estimators when the data fidelity term is Gaussian and the prior distribution is log-concave. Under some assumptions, the posterior mean estimator~\eqref{eqt: pm_estimator} and minimum mean squared error in~\eqref{eqt: MSE_def} associated to the posterior distribution~\eqref{eqt: post_dist} can be computed using the spatial gradient and Laplacian of the solution to the viscous HJ PDE~\eqref{eqt: 2ndorder_HJ} via the representation formulas~\eqref{eqt: prop1_relation_bayes_1} and~\eqref{eqt: prop1_relation_bayes_2}, respectively. We also proposed an analogue of the min-plus algebra technique designed for certain first-order HJ PDEs tailored to viscous HJ PDEs that enable us to compute posterior mean estimators with Gaussian fidelity term and prior that involves the sum of $m$ log-concave priors, i.e., to certain mixture models. The corresponding posterior mean estimator with non-convex regularization $J$ of the form of~\eqref{eqt: smoothminJi} can then be computed using the representation formulas~\eqref{eqt: i_pm_estimator} and posterior mean estimators~\eqref{eqt: pm_estimator} with convex regularization terms $J_i$.

Let us emphasize again that the proposed min-plus algebra technique for computations directly applies only for moderate $m$ in \eqref{eqt: minJi}. It would be of great interest to identify classes of non-convex regularizations for which novel numerical algorithms based on the min-plus algebra technique would not require to compute solutions to all $m$ convex sub-problems. To our knowledge, there is no available result in the literature on this matter.